\documentclass[a4paper,12pt]{article}
\usepackage{amssymb,amsfonts,amsmath,amsthm}
\usepackage{latexsym}
\usepackage{epsfig}
\newcommand{\RR}{{\mathbb{R}}}

\newcommand{\CC}{{\mathbb{C}}}

\newcommand{\NN}{{\mathbb{N}}}

\title{On two-dimensional shape-preserving approximation}
\author{Danylo Radchenko}
\newtheorem{prop}{Proposition}
\newtheorem{thm}{Theorem}
\newtheorem{prob}{Problem}
\newtheorem{prope}{Property}
\newtheorem{lem}{Lemma}

\begin{document}
\maketitle

\smallskip

\begin{abstract}
In this paper we investigate a problem of approximation of continuous mappings by smooth mappings with nonnegative Jacobian.\end{abstract}
%=======================================================================

\section{Introduction and main results} \label{sec:intro}
The following problem was communicated to me by Professor I.A. Shevchuk, who in turn was asked about it by Professor V.G. Krotov:
\begin{prob}
 What are the necessary and sufficient conditions for a continuous mapping $f:[0,1]^2 \rightarrow \RR^2$ to be uniformly approximable by $C^1$-smooth mappings with nonnegative Jacobian?
\end{prob}

In this paper we will provide some necessary and some sufficient conditions for this problem, as well as for the related problem.
\begin{prob}
 What are the necessary and sufficient conditions for a continuous mapping $f:[0,1]^2 \rightarrow \RR^2$ to be uniformly approximable by $C^1$-smooth mappings with positive Jacobian?
\end{prob}

We begin with some notation. A connected open subset of $\RR^2$ is called a \textit{domain}. We will only consider domains with piecewise smooth boundary. Denote by $|x|$ the Euclidean norm of $x\in \RR^2$. For a continuous mapping $f:K\rightarrow \RR^2$, where $K\subset \RR^2$ is a compact set, define $\left\|f\right\|_{K} = \max_{x\in K}|f(x)|$. Whenever it is clear what $K$ is, we write $\left\|f\right\|$ instead of $\|f\|_{K}$. For $x\in \RR^2$ and $r>0$ define $B_r(x)$ to be an open disk of radius $r$ centered at $x$. For any $A\subset \RR^2$ define the open neighborhood $O_r(A)=\bigcup_{x\in A}B_r(x)$. For sets $F_1,F_2\subset \RR^2$, by definition put $d(F_1,F_2)=\inf_{x\in F_1,y\in F_2}|x-y|$. For a set $F\subset \RR^2$ define $diam(F)=\sup_{x,y\in F}|x-y|$. We will identify $\RR^2$ with $\CC$ whenever it is convenient. For arbitrary set $A\subset \RR^2$, smooth (analytic, complex analytic) mappings on $A$ are considered to be defined on some open set $U\supset A$.

We will also consider both problems for mappings defined on $\overline{\Omega}$, where $\Omega$ is a bounded domain.

We will give the necessary conditions for both problems in terms of the Brouwer degree. For a nice overview of the Brouwer degree theory see \cite{Brouwer}.

Let $U$ be a bounded open subset of $\RR^N$ and $f:\overline{U} \rightarrow \RR^N$ be a continuous mapping. Then for any point $p\in \RR^N\setminus f(\partial U)$ one can define an integer $deg(f,U,p)$, which has the following properties (see \cite{Brouwer}, Th.1.2.6):

(i) If $f \in C^1(\overline{U})$ and $p$ is a regular value (that is for every $x\in  f^{-1}(p)$ we have $J_f(x)\neq 0$) for $f$ then
$$deg(f,U,p)=\sum_{x\in f^{-1}(p)}sgn J_f(x);$$

(ii) If $g:\overline{U} \rightarrow \RR^N$ is another mapping and $\left\|g-f\right\|_{\partial U}<d(p,f(\partial U))/2$, then $deg(f,U,p)=deg(g,U,p)$;

(iii) For a fixed mapping $f$, the value $deg(f,U,p)$ depends only on the connected component of $\RR^N\setminus f(\partial U)$ the point $p$ belongs to, and $deg(f,U,p)=0$ if $p\in \RR^N\setminus f(\overline{U})$;

(iv) If $\partial U$ is a simple closed curve, then $deg(f,U,p)$ is the total number of times that curve $f(\partial U)$ travels counterclockwise around the point $p$ (winding number).

From these properties we see that if the mapping $f:\overline{\Omega}\rightarrow \RR^N$ can be approximated by $C^1$-smooth mappings with nonnegative Jacobian, then $f$ must satisfy the following property.

\begin{prope}
For every open subset $U\subset\Omega$ and any $p\in f(U)\setminus f(\partial U)$ we have $deg(f,U,p)\geq 0$.
\label{prop:cond}
\end{prope}

From now we will concentrate on dimension $N=2$, as for $N\geq 3$ the situation is much more complicated.

Recall that a continuous mapping $f:X\rightarrow Y$ is called \textit{light} if $f^{-1}(y)$ is totally disconnected (connected components in $f^{-1}(y)$ are the one-point sets) for all $y\in Y$. For light mappings it turns out that if we require strict inequality in Property 1, then it becomes a sufficient condition for Problem 1.

\begin{thm}
Let $\Delta$ be a bounded domain and $f:\Delta\rightarrow \RR^2$ be a light continuous mapping. Suppose that for every open set $U$ with $\overline{U}\subset\Delta$ and any $p\in f(U)\setminus f(\partial U)$ we have $deg(f,U,p)>0$.

Then for any domain $\Omega$ such that $\overline{\Omega}\subset\Delta$ and each $\varepsilon > 0$, there exists a $C^\infty$-smooth mapping $g:\overline{\Omega}\rightarrow \RR^2$ with nonnegative Jacobian such that $\left\|f-g\right\|_{\overline{\Omega}}<\varepsilon$.
\label{thm:main2}\end{thm}

For the Problem 2 there is a very simple sufficient condition.
\begin{thm}
Suppose that a continuous mapping $f:\overline{\Omega}\rightarrow \RR^2$ is locally one-to-one. Then for each $\varepsilon > 0$ there exists a polynomial mapping $p:\overline{\Omega}\rightarrow \RR^2$ with nonzero Jacobian such that $\left\|f-p\right\|<\varepsilon$.
\label{thm:main}
\end{thm}

In an old paper \cite{FranklinWiener} there is a result about approximation of homeomorphisms similar to Theorem \ref{thm:main}. We only require the mapping to be locally univalent, and our approach can be used to establish the result of \cite{FranklinWiener}.

In the next section we will prove Theorem \ref{thm:main2} and Theorem \ref{thm:main}. The proof of Theorem \ref{thm:main} is rather technical and relies on piecewise linear approximation. In contrast, the proof of Theorem \ref{thm:main2} is more delicate and is based on the classification of light open mappings.

\section{Proof of main results} \label{sec:smooth}

Recall some definitions from piecewise linear topology. A simplicial complex $\mathcal{K}$ is a finite set of closed triangles in $\RR^2$ such that the intersection of any two triangles $\sigma_1, \sigma_2 \in \mathcal{K}$ is their common edge, common vertex, or an empty set. For a simplicial complex $\mathcal{K}$ denote by $|\mathcal{K}|$ the union of its triangles. A simplicial complex $\mathcal{L}$ is called a subdivision of a simplicial complex $\mathcal{K}$ if every triangle of $\mathcal{L}$ is contained in some triangle of $\mathcal{K}$ and $|\mathcal{K}|=|\mathcal{L}|$.

A mapping $f:|\mathcal{K}|\rightarrow \RR^2$ is called \emph{linear} if it is linear on every simplex of $\mathcal{K}$. A mapping $f:|\mathcal{K}|\rightarrow \RR^2$ is called \emph{piecewise linear} if it is linear for some subdivision of $\mathcal{K}$.

We will need the following result on piecewise linear approximation.
\begin{lem}
Let $\mathcal{K}$ be a simplicial complex and $f:|\mathcal{K}|\rightarrow \RR^2$ be a locally one-to-one continuous mapping. Then for each $\varepsilon > 0$ there exists a locally one-to-one piecewise linear mapping $h:|\mathcal{K}|\rightarrow \RR^2$ such that $\left\|f-h\right\|<\varepsilon$.
\label{lem:pllu}
\end{lem}

The proof is very similar to the proof of Theorem 6.3 from \cite{GeometricTopology}, but for the sake of completeness we give it in the Appendix.

In the following lemma denote by $I(C)$ the interior of a simple closed curve $C$ (given by the Jordan curve theorem).

\begin{lem}
Let $C$ be a unit circle, $r>0$ and $\phi:O_{r}(C)\rightarrow \RR^2$ be a $C^1$-smooth diffeomorphism. Suppose that $\phi(I(C)\cap O_{r}(C))\subset I(\phi(C))$ and $0\in I(\phi(C))\setminus \phi(O_{r}(C)).$ Then for some $\varepsilon\in (0,r)$ there exists a $C^1$-smooth diffeomorphism $\Phi:B_1(0)\cup O_{\varepsilon}(C)\rightarrow\RR^2$ such that $\phi|_{O_{\varepsilon}(C)}=\Phi|_{O_{\varepsilon}(C)}$.
\label{lem:sm}
\end{lem}
\textit{Proof:} The proof is found in \cite{HM2} and \cite{HM1}. \qed

\begin{lem}
 Let $A,B:\RR^2\rightarrow\RR^2$ be linear maps. Suppose that $\det(A)>0$, $\det(B)>0$, and there exists a nonzero vector $x$ such that $Ax=Bx$. Then for any $\alpha>0$, $\beta>0$ we have $\det(\alpha A+\beta B)>0$ .
\label{lem:la}
\end{lem}

\textit{Proof:} The proof is straightforward. \qed

\begin{lem}
Let $\mathcal{K}$ be a simplicial complex, $f:|\mathcal{K}|\rightarrow \RR^2$ be a locally one-to-one piecewise linear mapping. Then for each $\varepsilon>0$ there exists a $C^1$-smooth mapping $g:|\mathcal{K}|\rightarrow \RR^2$ with nonzero Jacobian such that $\left\|f-g\right\|<\varepsilon$.
\label{lem:smooth}
\end{lem}

\textit{Proof:} It is easily proved that we can extend the map $f$ to some simplicial complex $\mathcal{K}'\supset \mathcal{K}$, so that $|\mathcal{K}|$ is contained in the interior of $|\mathcal{K}'|$ and the extension of $f$ is locally one-to-one.

A standard nonnegative $C^{\infty}$-function with compact support is given by:

$$ \omega(x)=\begin{cases}  ce^{-1/(1-|x|^2)},& \mbox{if } |x|<1 \\ 0,& \mbox{else} \end{cases}$$

$$ \omega_{\delta}(x)=\delta^{-2}\omega(x/\delta) .$$

Note that $\omega$ is normalized, that is $\int_{\RR^2}\omega(x)=1$.

The convolution $f \ast \omega_{\delta}$ is a $C^{\infty}$-smooth mapping that uniformly converges to $f$, as $\delta\rightarrow 0$. It is obvious that if $f$ is linear on $B_{\delta}(x)$, then $f \ast \omega_{\delta}(x)=f(x)$. Note that if $B_{\delta}(x)$ only intersects two triangles of $\mathcal{K}$, then $D(f\ast \omega_{\delta})(x)$ is a convex combination of two linear maps, and it is easy to see that these maps satisfy conditions of Lemma \ref{lem:la}. Therefore, in this case $J_{f \ast \omega_{\delta}}(x)>0$.

Let $\delta_1>0$ be such that for every vertex $v_i$ of $\mathcal{K}$ we have $diam(f(B_{\delta_1}(v_i)))<\varepsilon$, $B_{\delta_1}(v_i)$ are pairwise disjoint, and $B_{\delta_1}(v_i)$ intersects $\mathcal{K}^1$ (the set of edges of all triangles in $\mathcal{K}$) only at edges that have $v_i$ as an endpoint. Let $\delta$ be such that $|f(x)-f(y)|<\varepsilon$ whenever $|x-y|<\delta$.  Let $\delta_2<\delta$ be such that for $x\in \mathcal{K}\setminus \bigcup_{v_i \in K^0} B_{\delta_1/2}(v_i)$ we have that $B_{\delta_2}(x)$ intersects $\mathcal{K}$ at not more than two triangles. It is easy to check that $f \ast \omega_{\delta_2}$ is one-to-one on some neighborhood of $\partial B_{2\delta_1/3}(v_i)$. Therefore, using Lemma \ref{lem:sm} we complete the construction by extending $g=f \ast \omega_{\delta_2}$ from $\mathcal{K}\setminus \bigcup_{v_i \in \mathcal{K}^0} B_{2\delta_1/3}(v_i)$ to the whole $\mathcal{K}$. \qed
\\

Now we prove Theorem \ref{thm:main2} and Theorem \ref{thm:main}. To avoid technical details, we only prove Theorem \ref{thm:main} in case when $\overline{\Omega}$ is equal to $|\mathcal{K}|$ for some simplicial complex $\mathcal{K}$.

\textit{Proof of Theorem \ref{thm:main}:} From Lemma \ref{lem:pllu} and Lemma \ref{lem:smooth} we see that there exists a $C^1$-smooth mapping $g:\overline{\Omega}\rightarrow \RR^2$ with nonzero Jacobian such that $\|f-g\|<\varepsilon/2$. Since $\overline{\Omega}$ is a compact set, there exists $\delta>0$ such that $|J_g(x)|>\delta$ for all $x\in \overline{\Omega}$. Define $M=\max\{\|\frac{\partial g_i}{\partial x_j}\|: 1\leq i,j \leq 2\}$. By a theorem on simultaneous approximation(see \cite{SimApprox}, Th.1) there exists a polynomial map $p:\overline{\Omega}\rightarrow \RR^2$ such that $\|g-p\|<\varepsilon/2$ and $\|\frac{\partial g_i}{\partial x_j}-\frac{\partial p_i}{\partial x_j}\|< \max(2M,\frac{\delta}{8M})$. Then it is easy to see that $\|f-p\|<\varepsilon$ and $|J_p(x)|>0$. \qed

\textit{Proof of Theorem \ref{thm:main2}:} Recall that a mapping $F:X\rightarrow \RR^2$ is called \textit{quasi-open} if for any $y\in F(X)$ and any open set $V$ containing a compact component of $F^{-1}(y)$, $y$ is interior to $F(V)$. Note that our mapping $f$ is quasi-open. Indeed, for every $x \in \Delta$ and $V$ containing a compact component of $f^{-1}(f(x))$, there is some open subset $V_0\subset V$ such that $\partial V_0 \cap f^{-1}(f(x))=\emptyset$. Then we have that $deg(f,V_0,f(x))>0$. Therefore, by property (iii) of the Brouwer degree we have that $f(V_0)$ contains some neighborhood of $f(x)$. Hence, $f$ is quasi-open. Since $f$ is also light, we have (see \cite{Whyburn}, pp.110-113) that $f$ is open. Since $f$ is open and light, by theorem of Stoilow (\cite{Whyburn}, p. 103) $f$ is topologically equivalent to a complex analytic mapping $h:U_1\rightarrow U_2\subset \RR^2$. Then there are homeomorphisms $s_1:\Delta\rightarrow U_1$ and $s_2:U_2\rightarrow s_2(U_2)\subset \RR^2$ such that $f=s_2hs_1$. By Theorem \ref{thm:main} there are polynomial mappings $p_1,p_2$ with nonzero Jacobian such that $\|s_2hs_1-s_2hp_1\|_{\overline{\Omega}}<\varepsilon/2$ and $\|s_2-p_2\|_{h(p_1(\overline{\Omega}))}<\varepsilon/2$. Then $\|f-p_2hp_1\|_{\overline{\Omega}}<\varepsilon$. Since $f$ preserves orientation, $p_1$ and $p_2$ must have Jacobian of the same sign. Therefore, $p_2hp_1$ is a $C^\infty$-smooth mapping with nonnegative Jacobian. \qed

\section{Negative result} \label{sec:negres}

Here we prove a negative result that shows the difference between problems 1 and 2. Let the mapping $f:\overline{B_1(0)}\rightarrow \CC$ be given by $f(z)=z^2$.
\begin{thm}
Let the mapping $f$ be defined as above. Then for any $C^1$-smooth mapping $g:\overline{B_1(0)}\rightarrow \RR^2$ with strictly positive Jacobian we have $\left\|f-g\right\|\geq 1/4$.
\label{thm:neg1}
\end{thm}

\textit{Proof:} Suppose that there exists a mapping $g$ such that $\left\|f-g\right\|<1/4$. Then for every $x \in \overline{B_{1/2}(0)}$ we have $deg(g,B_1(0),x)=2$. Since g has positive Jacobian, there are exactly two different solutions $y_1,y_2$ to $g(y)=x$, with $x$ in $\overline{B_{1/2}(0)}$.

Consider the set $U=g^{-1}(B_{1/2}(0))$. We see that $g:U\rightarrow B_{1/2}(0)$ is a covering map and $U$ is a double cover of $B_{1/2}(0)$. Since $B_{1/2}(0)$ is simply connected, we get by classification of covering spaces (see \cite{Hatcher}, Th.1.38) that $U$ is a disjoint union of two homeomorphic copies of $B_{1/2}(0)$. Then $U=U_1 \cup U_2$ and $g$ maps $U_i$ homeomorphically onto $B_{1/2}(0)$, so we can define two inverse maps $y_1,y_2$ from $B_{1/2}(0)$ to $U$.

From $\left\|f-g\right\|<1/4$ we have $|y_i^2(z)-z|< 1/4$. Consider the map $\gamma(\phi)=e^{i\phi}/(2+\varepsilon)$. Then the points $y_i(\gamma(\phi))$ are contained in a disjoint union of disks $B_{1/2}(e^{i\phi/2}/\sqrt{2+\varepsilon}))$ and $B_{1/2}(e^{i(\phi/2+\pi)}/\sqrt{2+\varepsilon}))$. It is easy to verify that each disk contains exactly one of points $y_1(\gamma(\phi)),y_2(\gamma(\phi))$. Then as we continuously change $\phi$ from $0$ to $2\pi$, $y_1(\gamma(\phi))$ ends up in a different disk, so we must have $y_1(\gamma(0))=y_2(\gamma(0))$. This contradiction concludes the proof.  \qed
\\

Since $J_f(x_1,x_2)=x_1^2+x_2^2 \geq 0$, we have a polynomial mapping with nonnegative Jacobian that cannot be approximated by $C^1$-smooth mappings with positive Jacobian.

\section{Appendix} \label{sec:PLap}
Here we prove Lemma 1. We need the following result.

\begin{prop}(PL Schoenflies theorem)
Let $\Delta$ be a nondegenerate triangle and $f:\partial \Delta \rightarrow \RR^2$ be a piecewise linear (PL) homeomorphism. Then there exists a PL homeomorphism $h:\Delta \rightarrow \RR^2$ such that $h|_{\partial \Delta} =f$.
\end{prop}
The proof is found in \cite{GeometricTopology}.

Combinatorial distance between vertices of a simplicial 1-complex is a minimal number of edges in any path connecting them. Denote by $dist(v_i,v_j)$ the combinatorial distance between $v_i$ and $v_j$. For sets of vertices $U,V$ define $dist(U,V):= \max_{u\in U,v\in V}dist(u,v)$. Combinatorial diameter of $U$ is then defined as $dist(U,U)$.

\begin{lem}
Let $\mathcal{K}^1$ be a simplicial 1-complex, $d\in \NN, d\geq 2$. Suppose that a continuous mapping $f:|\mathcal{K}^1|\rightarrow \RR^2$ is one-to-one on every subcomplex of $\mathcal{K}^1$ with combinatorial diameter $\leq d$. Then for each $\varepsilon > 0$ there exists a piecewise linear mapping $g:|\mathcal{K}^1|\rightarrow \RR^2$ such that $g$ is one-to-one on every subcomplex of $\mathcal{K}^1$ with combinatorial diameter $\leq d$, $g$ satisfies $\left\|f-g\right\|<\varepsilon$, and $f(v)=g(v)$ for each vertex $v$ of $\mathcal{K}^1$. \label{lem:pl1}
\end{lem}
\textit{Proof:} We may assume that $diam(f(v_iv_j))<\varepsilon/3$ for any edge $v_iv_j$ of $\mathcal{K}^1$ (f is uniformly continuous and by taking small subdivisions of $\mathcal{K}^1$ with equal number of parts on each edge, we can ensure that f still has the property for a multiple of $d$). Let $v_i$ be vertices of $\mathcal{K}^1$, $w_i=f(v_i)$ and $A_{ij}=f(v_iv_j)$. Then we have $diam(A_{ij})<\varepsilon/3$, so for all $x,y \in O_{\varepsilon/3}(A_{ij})$ we have $d(x,y)<\varepsilon$. Define $N_i = O_{\varepsilon_i}(w_i)$, where $\varepsilon_i$ are small enough, so that:

(i) $\overline{N_i} \cap \overline{N_j}=\emptyset$ whenever $dist(v_i,v_j)\leq d$;

(ii) $\varepsilon_i < \varepsilon/3$;

(iii) For any three vertices $v_i,v_k,v_j$ such that $dist(v_i,v_j)\leq d$, $dist(v_i,v_k)\leq d$ and $v_jv_k$ is an edge, we have that $\overline{N_i} \cap A_{kj} \neq \emptyset$ exactly when $v_k=v_i$ or $v_j=v_i$.

Let $x_{ij}$ be the last point of $A_{ij}$ (in the order from $w_i$) that lies in $\overline{N_i}$. Let $x_{ij}'$ be the first point of $A_{ij}$ that follows $x_{ij}$ and belongs to $\overline{N_j}$. Let $A_{ij}'$ be the arc from $x_{ij}$ to $x_{ij}'$ in  $A_{ij}$. Then arcs $A_{ij}'$ and $A_{kl}'$ are disjoint whenever $dist(v_iv_j,v_kv_l)\leq d$. Next, take $\delta$-neighborhoods of $A_{ij}'$ with $\delta<\varepsilon/3$, so that they are disjoint whenever corresponding $A_{ij}'$ are disjoint. Then for each edge $v_{ij}$ there is a broken line $B_{ij}$ in $O_{\delta}(A_{ij}')$ that joins $x_{ij}$ and $x_{ij}'$ (cf. \cite{GeometricTopology}, Th.6.1). Therefore, $B_{ij}$ and $B_{kl}$ are disjoint whenever $dist(v_iv_j,v_kv_l)\leq d$.

Let $y_{ij}$ be the last point of $B_{ij}$ that lies in $\overline{N_i}$, and $y_{ij}'$ be the first point of $B_{ij}$ that follows $y_{ij}$ and belongs to $\overline{N_j}$. Now define $B_{ij}'$ to be a part of broken line $B_{ij}$ from $y_{ij}$ to $y_{ij}'$. Finally, let $B_{ij}''=w_iy_{ij} \cup B_{ij}'\cup y_{ij}'w_j$. The broken line $B_{ij}''$ connects $w_i$ to $w_j$. The edges $v_iv_j$ and $v_kv_l$ can only intersect at endpoints whenever $dist(v_iv_j,v_kv_l)\leq d$. We also have that $B_{ij}'' \subset O_{\varepsilon/3}(A_{ij})$. We construct $g:|\mathcal{K}^1|\rightarrow \RR^2$ by defining each mapping $g|_{v_iv_j}$ to be a piecewise linear homeomorphism that sends $v_iv_j$ to $B_{ij}''$, $v_i$ to $w_i$, and $v_j$ to $w_j$. Then $g$ has the desired property and interpolates $f$ at vertices of $\mathcal{K}^1$. To show that $g$ is an $\varepsilon$-approximation for $f$, notice that for $x \in v_iv_j$ both $f(x)$ and $g(x)$ lie in $O_{\varepsilon/3}(A_{ij})$, so that $|f(x)-g(x)|<\varepsilon$. $\qed$
\\

\textit{Proof of Lemma 1:} First, for every point $x \in |\mathcal{K}|$ consider neighborhood $O_{\varepsilon}(x)$ such that $f|_{O_{\varepsilon}(x)}$ is one-to-one. The family of these neighborhoods covers $|\mathcal{K}|$, so we may choose a finite subcover. Now let $\delta$ be the Lebesgue number of this cover. Let $\mathcal{L}$ be a subdivision of $\mathcal{K}$ such that for every triangle $\sigma \in \mathcal{L}$, we have $diam(\sigma)<\delta/3$, and $diam(f(\sigma))<\varepsilon/3$. For every triangle $\sigma \in \mathcal{L}$ let $D(\sigma)$ be set of vertices of $\mathcal{L}$ at combinatorial distance not more than 3 from $\sigma$, excluding the vertices of $\sigma$. Next, let $\theta_{\sigma}=\min\left\{\varepsilon/3,d(f(\sigma),f(D(\sigma)))\right\}$. Finally, let $\delta_1=\min_{\sigma}{\theta_{\sigma}}$.

Now take $g_1:|\mathcal{L}^1| \rightarrow \RR^2$ to be a $\delta_1$-approximation to $f$ from Lemma \ref{lem:pl1} with $d=3$. Next, using PL Schoenflies theorem, extend this function to $g:|\mathcal{L}|\rightarrow \RR^2$ in such a way that it will be one-to-one on each triangle. First, we prove that $g$ is an $\varepsilon$-approximation for $f$. For every triangle $\sigma$ in $\mathcal{L}$ we know, that $diam(f(\sigma))<\varepsilon/3$. Since $\delta_1<\varepsilon/3$, we have that $g(\partial\sigma) \subset O_{\varepsilon/3}(f(\sigma))$. Since $g|_{\sigma}$ is a homeomorphism, we also have that $g(\sigma) \subset O_{\varepsilon/3}(f(\sigma))$. Thus, for any $x\in \sigma$, $f(x)$ and $g(x)$ both lie in $O_{\varepsilon/3}(f(\sigma))$, so $|f(x)-g(x)|<\varepsilon$.

To prove that this PL mapping is locally univalent, consider the subdivision $\mathcal{L}_1$ of $\mathcal{L}$, on which the mapping is linear. Let $\sigma_1'$ and $\sigma_2'$ be two different triangles of $\mathcal{L}_1$ which lie in triangles $\sigma_1$ and $\sigma_2$ of $\mathcal{L}$ correspondingly, and have non-empty intersection. If $\sigma_1=\sigma_2$, there is nothing to prove. If $\sigma_1'$ and $\sigma_2'$ have a common edge, then $\sigma_1$ and $\sigma_2$ also have a common edge. If $\sigma_1'$ and $\sigma_2'$ have a common vertex, then $\sigma_1$ and $\sigma_2$ also have a common vertex. So, combinatorial diameter of $\sigma_1 \cup \sigma_2$ is 2, and so $g_1$ is one-to-one on $\partial\sigma_1 \cup \partial\sigma_2$. So $g(\partial\sigma_1)\cap g(\partial\sigma_2)=g(\partial\sigma_1 \cap \partial\sigma_2)$, and this means that either $g(\sigma_1)\subset g(\sigma_2)$ or $g(\sigma_1)\cap g(\sigma_2)=g(\sigma_1 \cap \sigma_2)$. But if we have $g(\sigma_1)\subset g(\sigma_2)$, then one of the vertices if $\sigma_2$ maps inside $g(\sigma_1)$. But $g(\sigma_1) \subset O_{\theta_{\sigma_1}}(f(\sigma_1))$, so this contradicts to definition of $\theta_{\sigma_1} \leq d(f(\sigma_1),f(D(\sigma_1)))$, because the other vertex of $\sigma_2$ is in $D(\sigma_1)$. $\qed$

\end{document}